\documentclass{amsart}
\usepackage{amsmath,amsthm,amsfonts,amssymb}
\def\({\left(}
\def\){\right)}

\def\cb{{\mathcal B}}

\theoremstyle{plain}   
\begingroup 
\newtheorem*{thm*}{Theorem}          
\newtheorem*{cor*}{Corollary}        
\endgroup

\theoremstyle{definition}

\newtheorem*{rem*}{Remark}

\theoremstyle{remark}
\theoremstyle{definition}

\numberwithin{equation}{section}

\begin{document}


\title[Proofs of the Ergodic Theorem and Maximal Ergodic Theorem]
{Easy and nearly simultaneous proofs \\of the Ergodic
Theorem and Maximal Ergodic Theorem}

\author{Karl Petersen}

\address{Department of Mathematics,
CB 3250, Phillips Hall,
         University of North Carolina,
Chapel Hill, NC 27599 USA}
\email{petersen@math.unc.edu}


\date{\today}


\maketitle


This note arose from a conversation I had with Mike Keane after his
 talk at E.S.I. in Vienna on a simple proof of the Ergodic
Theorem. 
Thanks to E. Lesigne and X. Mela for inducing several clarifications.

Let $(X, \cb ,\mu )$ be a probability space, $T:X \to X$ a (possibly
noninvertible) measure-preserving transformation, and $f \in L^1(X, \cb
, \mu )$. Let
\[
A_kf=\frac{1}{k}\sum_{j=0}^{k-1}fT^j, \quad f_N^*= \sup_{1\leq k \leq
N} A_kf, \quad f^*=\sup_N f_N^*, \quad \text{ and } \overline A =
\limsup _{k \to \infty} A_kf .
\]
When $\lambda$ is a constant, the following result is the Maximal
Ergodic Theorem. Choosing $\lambda = \overline A - \epsilon$ covers
most of the proof of the Ergodic Theorem.

\begin{thm*} Let $\lambda$ be an invariant ($\lambda \circ T=\lambda$ a.e.) 
function on $X$ with $\lambda ^+ \in L^1$. Then
\begin{equation*}
\int_{\{f^*>\lambda\}}(f-\lambda) \geq 0 .
\end{equation*}
\end{thm*}

\begin{proof} 
We may assume that $\lambda \in L^1\{ f^*>\lambda \}$, since
otherwise
\[
\int_{ \{ f^*> \lambda \} }(f-\lambda ) = \infty \geq 0. 
\]
But then actually $\lambda \in L^1(X)$, since on $\{ f^* \leq
\lambda\}$ we have $f \leq \lambda$, so that on this set $\lambda^-
\leq -f + \lambda^+$, which is integrable.

Assume first that $f \in L^\infty$. Fix $N=1,2,\dots$, and let
\begin{equation*}
E_N=\{f_N^*>\lambda \} .
\end{equation*}
Notice that
\begin{equation*}
(f-\lambda ) \chi_{E_N} \geq (f-\lambda) ,
\end{equation*} 
since $x \notin E_N$ implies $(f - \lambda )(x) \leq 0$. 
Thus for a very large $m \gg N$, we can break up
\begin{equation*}
\sum_{k=0}^{m-1}(f-\lambda  )\chi_{E_N}(T^kx)
\end{equation*}
into convenient strings of terms as follows. There is maybe an initial
string of $0$'s during which $T^kx \notin E_N$. Then there is a
first time $k$ when $T^kx \in E_N$, which initiates a string of no
more than $N$ terms, the sum of which is positive (using on each of
these terms the fact that $(f-\lambda ) \chi_{E_N} \geq (f-\lambda)$).
Beginning after the last term in this string, we repeat the previous
analysis, finding maybe some $0$'s until again some $T^kx \in E_N$
initiates another string of no more than $N$ terms and with positive
sum. The full sum of $m$ terms may end in the middle of either of
these two kinds of strings ($0$'s, or having positive sum). Thus we
can find $j=m-N+1,\dots, m$ such that 
\begin{equation*}
\sum_{k=0}^{m-1}(f-\lambda )\chi_{E_N}(T^kx) \geq
\sum_{k=j}^{m-1}(f-\lambda )\chi_{E_N}(T^kx) \geq -N(\| f \|_\infty +
\lambda ^+(x)) .
\end{equation*}
Integrating both sides, dividing by $m$, and letting $m \to \infty$ gives
\begin{gather*}
m \int_{E_N}(f - \lambda) \geq  -N(\| f \|_\infty +
\| \lambda ^+\|_1) , \\
\int_{E_N}(f-\lambda ) \geq \frac{-N}{m}(\| f \|_\infty +
\| \lambda ^+\|_1) , \\
\int_{E_N}(f-\lambda ) \geq 0.
\end{gather*}
Letting $N \to \infty$ and using the Dominated Convergence Theorem
concludes the proof for the case $f \in L^\infty$.

To extend to the case $f \in L^1$, for $s=1,2,\dots$ let $\phi_s = f
\cdot \chi_{\{ |f| \leq s\}}$, so that $\phi_s \in L^\infty$ 
and $\phi_s \to f$ a.e. and in $L^1$.
Then for fixed
$N$
\begin{equation*}
(\phi_s)_N^* \to f_N^* \quad \text{a.e. and in } L^1 \quad\text{ and }\quad
\mu \left(\{(\phi_s)_N^* >\lambda \} \bigtriangleup \{ f_N^* >
  \lambda\} \right) \to 0.
\end{equation*}
Therefore
\[
0 \leq \int_{\{(\phi_s)_N^* > \lambda\} } (\phi_s - \lambda) \to \int
_{\{ f_N^* > \lambda \} } (f -\lambda) ,
\]
again by the Dominated Convergence Theorem. The full result follows by
letting $N \to \infty$.
\end{proof}

\begin{cor*}[Ergodic Theorem] 
The sequence $(A_kf)$ converges a.e..
\end{cor*} 
\begin{proof}
It is enough to show that
\begin{equation*}
\int \overline A \leq \int f  .
\end{equation*}
For then, letting $\underline A = \liminf A_kf$, applying this to $-f$ gives
\begin{equation*}
-\int \underline A \leq - \int f  ,
\end{equation*}
so that
\begin{equation*}
\int \overline A \leq \int f \leq \int \underline A \leq \int \overline A ,
\end{equation*}
and hence
\begin{equation*} 
\int( \overline A - \underline A) = 0 , \quad\text{ so that } 
\overline A = \underline A \text{  a.e..}
\end{equation*}

Consider first $f^+$ and its associated $\overline A$, denoted by
$\overline A(f^+)$.
For any invariant function $\lambda < \overline
  A(f^+)$ 
such that
  $\lambda^+  \in L^1$, for example $\lambda = \overline A(f^+) \wedge n
   - 1/n$, we have  
$\{(f^+)^* > \lambda\} = X$, so the Theorem gives
\begin{equation*}
\int f^+ \geq \int \lambda \nearrow \int \overline A(f^+) .
\end{equation*}
Thus $(\overline A)^+ \leq
\overline A(f^+)$ is integrable (and, by a similar argument, so is  
$(\overline A)^- \leq \overline A(f^-)$.)

Now let $\epsilon > 0$ be arbitrary and apply the Theorem to 
 $\lambda = \overline A - \epsilon$ to conclude that 
\begin{equation*}
\int f \geq \int \lambda \nearrow \int \overline A .
\end{equation*}
\end{proof}

\begin{rem*}
Keane's proof was a development of the Katznelson-Weiss proof
\cite{KW} based on Kamae's nonstandard-analysis proof \cite{Kamae} and
presented also in the Bedford-Keane-Series collection \cite{BKS}. That
proof is given in a paper by Keane \cite{Keane} and has been extended
to deal also with the Hopf Ratio Ergodic Theorem and with the case of
higher-dimensional actions \cite{KK}. The proof given above yields
both the Pointwise and Maximal Ergodic Theorems essentially
simultaneously without adding any real complications. Roughly
contemporaneously with our conversation, Roland Zweim\" uller prepared
some preprints \cite{Z,ZZ} also giving short proofs based on the
Kamae-Katznelson-Weiss approach.  Without going too deep into the
complicated history of the Ergodic Theorem and Maximal Ergodic
Theorem, it is interesting to note some recurrences as the use of
maximal theorems arose and waned repeatedly.  After the original
proofs by von Neumann \cite{VN}, Birkhoff \cite{Birk}, and Khinchine
\cite{Kh}, the role and importance of the Maximal Lemma and Maximal
Theorem were brought out by Wiener \cite{W} and Yosida-Kakutani
\cite{YK}, making possible the exploration of connections with
harmonic functions and martingales.  Proofs by upcrossings followed an
analogous pattern. It also became of interest, for instance to allow
extension to new areas or new kinds of averages, again to prove the
Ergodic Theorem without resort to maximal lemmas or theorems, as in
the proof by Shields \cite{Shields} inspired by the Ornstein-Weiss
proof of the Shannon-McMillan-Breiman Theorem for actions of amenable
groups \cite{OW}, or in Bourgain's proofs by means of variational
inequalities \cite{Bourgain}.  Sometimes it was pointed out, for
example in the note by R. Jones \cite{RJ}, that these approaches could
also with very slight modification prove the Maximal Ergodic
Theorem. Of course there are the theorems of Stein \cite{Stein} and
Sawyer \cite{Sawyer} that make the connection explicit, just as the
transference techniques of Wiener \cite{W} and Calde\' ron \cite{CALD}
connect ergodic theorems with their analogues in analysis like the
Hardy-Littlewood Maximal Lemma \cite{HL}.  In many of the improvements
over the years, ideas and tricks already in the papers of Birkhoff,
Kolmogorov \cite{Kolm}, Wiener, and Yosida-Kakutani have continued to
play an essential role.
\end{rem*}

 \bibliographystyle{amsplain} \bibliography{ergthpf}
\end{document}